%
\catcode`@=11
%
%
\def\bibn@me{R\'ef\'erences}
\def\bibliographym@rk{\centerline{{\sc\bibn@me}}
	\sectionmark\section{\ignorespaces}{\unskip\bibn@me}
	\bigbreak\bgroup
	\ifx\ninepoint\undefined\relax\else\ninepoint\fi}
%
%
%
\let\refsp@ce=\ 
\let\bibleftm@rk=[
\let\bibrightm@rk=]
%
%
%
\def\numero{n\raise.82ex\hbox{$\fam0\scriptscriptstyle o$}~\ignorespaces}
%
%
\newcount\equationc@unt
\newcount\bibc@unt
\newif\ifref@changes\ref@changesfalse
\newif\ifpageref@changes\ref@changesfalse
\newif\ifbib@changes\bib@changesfalse
\newif\ifref@undefined\ref@undefinedfalse
\newif\ifpageref@undefined\ref@undefinedfalse
\newif\ifbib@undefined\bib@undefinedfalse
\newwrite\@auxout
%
%
\def\eqnum{\global\advance\equationc@unt by 1%
\edef\lastref{\number\equationc@unt}%
\eqno{(\lastref)}}
%
%
%
%
%
%
\def\re@dreferences#1#2{{%
	\re@dreferenceslist{#1}#2,\undefined\@@}}
\def\re@dreferenceslist#1#2,#3\@@{\def\next{#2}%
	\expandafter\ifx\csname#1@@\meaning\next\endcsname\relax
	??\immediate\write16
	{Warning, #1-reference "\next" on page \the\pageno\space
	is undefined.}%
	\global\csname#1@undefinedtrue\endcsname
	\else\csname#1@@\meaning\next\endcsname\fi
	\ifx#3\undefined\relax
	\else,\refsp@ce\re@dreferenceslist{#1}#3\@@\fi}
%
%
%
\def\newlabel#1#2{{\def\next{#1}\newl@bel#2}}
\def\newl@bel#1#2{%
	\expandafter\xdef\csname ref@@\meaning\next\endcsname{#1}%
	\expandafter\xdef\csname pageref@@\meaning\next\endcsname{#2}}
\def\label#1{{%
	\toks0={#1}\message{ref(\lastref) \the\toks0,}%
	\ignorespaces\immediate\write\@auxout%
	{\noexpand\newlabel{\the\toks0}{{\lastref}{\the\pageno}}}%
	\def\next{#1}%
	\expandafter\ifx\csname ref@@\meaning\next\endcsname\lastref%
	\else\global\ref@changestrue\fi%
	\newlabel{#1}{{\lastref}{\the\pageno}}}}
\def\ref#1{\re@dreferences{ref}{#1}}
\def\pageref#1{\re@dreferences{pageref}{#1}}
%
%
\def\bibcite#1#2{{\def\next{#1}%
	\expandafter\xdef\csname bib@@\meaning\next\endcsname{#2}}}
\def\cite#1{\bibleftm@rk\re@dreferences{bib}{#1}\bibrightm@rk}
%
%
\def\beginthebibliography#1{\bibliographym@rk
	\setbox0\hbox{\bibleftm@rk#1\bibrightm@rk\enspace}
	\parindent=\wd0
	\global\bibc@unt=0
	\def\bibitem##1{\global\advance\bibc@unt by 1
		\edef\lastref{\number\bibc@unt}
		{\toks0={##1}
		\message{bib[\lastref] \the\toks0,}%
		\immediate\write\@auxout
		{\noexpand\bibcite{\the\toks0}{\lastref}}}
		\def\next{##1}%
		\expandafter\ifx
		\csname bib@@\meaning\next\endcsname\lastref
		\else\global\bib@changestrue\fi%
		\bibcite{##1}{\lastref}
		\medbreak
		\item{\hfill\bibleftm@rk\lastref\bibrightm@rk}%
		}
	}
\def\endthebibliography{\egroup\par}
%
%
%
\def\@closeaux{\closeout\@auxout
	\ifref@changes\immediate\write16
	{Warning, changes in references.}\fi
	\ifpageref@changes\immediate\write16
	{Warning, changes in page references.}\fi
	\ifbib@changes\immediate\write16
	{Warning, changes in bibliography.}\fi
	\ifref@undefined\immediate\write16
	{Warning, references undefined.}\fi
	\ifpageref@undefined\immediate\write16
	{Warning, page references undefined.}\fi
	\ifbib@undefined\immediate\write16
	{Warning, citations undefined.}\fi}
%
%
\immediate\openin\@auxout=\jobname.aux
\ifeof\@auxout \immediate\write16
  {Creating file \jobname.aux}
\immediate\closein\@auxout
\immediate\openout\@auxout=\jobname.aux
\immediate\write\@auxout {\relax}%
\immediate\closeout\@auxout
\else\immediate\closein\@auxout\fi
%
%
\input\jobname.aux
\immediate\openout\@auxout=\jobname.aux
%
%

\def\bibn@me{R\'ef\'erences bibliographiques}
%
\def\bibliographym@rk{\bgroup}
%
%
\outer\def\bye{ 	\par\vfill\supereject\end}

\def\Q{{\bf {Q}}}

\def\Z{{\bf Z}} 
\def\R{{\bf R}}

\overfullrule=0pt

\magnification=1200

  \def\pro{\noindent {\bf{Proof : }}}

\def\house#1{\setbox1=\hbox{$\,#1\,$}
\dimen1=\ht1 \advance\dimen1 by 2pt \dimen2=\dp1 \advance\dimen2 by 2pt
\setbox1=\hbox{\vrule height\dimen1 depth\dimen2\box1\vrule}
\setbox1=\vbox{\hrule\box1}
\advance\dimen1 by .4pt \ht1=\dimen1
\advance\dimen2 by .4pt \dp1=\dimen2 \box1\relax}

\def\Card{{\rm Card}}

  \def\deg{{\rm deg}}

\def\disp{\displaystyle} 
  \def\noi{\noindent}

\def\build#1_#2^#3{\mathrel{\mathop{\kern 0pt#1}\limits_{#2}^{#3}}}

\def\date {le\ {\the\day}\ \ifcase\month\or 
janvier\or fevrier\or mars\or avril\or mai\or juin\or juillet\or
ao\^ut\or septembre\or octobre\or novembre\or 
d\'ecembre\fi\ {\oldstyle\the\year}}

\font\fivegoth=eufm5 \font\sevengoth=eufm7 \font\tengoth=eufm10

\newfam\gothfam \scriptscriptfont\gothfam=\fivegoth
\textfont\gothfam=\tengoth \scriptfont\gothfam=\sevengoth

\def\cqfd{\unskip\kern 6pt\penalty 500 \raise 0pt\hbox{\vrule\vbox 
to6pt{\hrule width 6pt \vfill\hrule}\vrule}\par}

\def\pro{\noindent {\it Proof. }}

\def\smallsquare{\vbox{\hrule\hbox{\vrule height 1 ex\kern 1 ex\vrule}\hrule}}
\def\cqfd{\hfill \smallsquare\vskip 3mm}

\def\Card{{\rm Card }\, }

\def\FF{{\bf F}}
\def\kk{{\bf k}}
\def\wh{\widehat}
\def\QI{{K^{(2)}}}
\def\Sp{{\rm Sp}}

\def\Hw{{\rm Hw}}
\def\v{{\rm v}}
\def\deg{{\rm deg} \, }

\centerline{}

\vskip 4mm

\centerline{
\bf Nonarchimedean quadratic Lagrange spectra }

\medskip

\centerline{\bf and continued fractions in power series fields}

\vskip 13mm

\centerline{
Y{\sevenrm ANN} B{\sevenrm UGEAUD}
\footnote{}{\rm
Lagrange spectrum, continued fractions, power series fields.}}

{\narrower\narrower
\vskip 15mm

\proclaim Abstract. {
Let $\FF_q$ be a finite field of order a positive power $q$ of a 
prime number. 
We study the nonarchimedean quadratic Lagrange spectrum defined
by Parkkonen and Paulin by considering the
approximation by elements of the orbit of a given quadratic
power series in $\FF_q((Y^{-1}))$, 
for the action by homographies and anti-homographies of
${\rm PGL}_2(\FF_q[Y])$ on $\FF_q((Y^{-1})) \cup \{\infty\}$. 
While their approach used geometric methods of group actions on Bruhat--Tits trees, 
ours is based on the theory of continued fractions in power series fields.
}

}

\vskip 12mm

\centerline{\bf 1. Introduction}

\vskip 6mm

For an irrational real number $\xi$, define $\lambda (\xi)$ in $(0, + \infty]$ by
$$
\lambda (\xi)^{-1} := \liminf_{p, q \in \Z, q \to + \infty} \, |q (q \alpha - p)|.
$$
The Lagrange spectrum ${\cal L}$ is the set of values taken by the function $\lambda$ 
at irrational real numbers. It is included in $[\sqrt{5}, + \infty]$ and has a rather complicated structure, 
which is not completely understood, despite some recent progress \cite{Gay17,Mo18}. 
The first values of ${\cal L}$ are $\sqrt{5}$ (sometimes called Hurwitz' constant) 
and $2 \sqrt{2}$, and its first accumulation point is $3$. 
In 1947 Hall \cite{Hall47} established that the whole interval $[6, + \infty]$ is contained in 
${\cal L}$. In 1975, Freiman \cite{Fr75} proved that the biggest half-line contained in ${\cal L}$ is 
$$
\biggl[ {2221564096 + 283748 \sqrt{462} \over 491993569}, + \infty \biggr].
$$
This half-line is called {\it Hall's ray}. The reader is directed to \cite{CuFl} for additional references (note 
also that, sometimes, the authors study the set of values 
taken by the function $1/ \lambda$). 

Analogous spectra have been defined and studied in various contexts, including 
inhomogeneous Diophantine approximation 
(see e.g. Cusick, Moran, and Pollington \cite{CuMoPo96} and Pinner and Wolczuk \cite{PiWo01}), 
Diophantine approximation in imaginary quadratic fields (Maucourant \cite{Mau03}), 
and in the setting of interval exchange transformations and the Teichm\"uller 
flow on moduli spaces of translation surfaces (Hubert, Marchese, and Ulcigrai \cite{HuMaUl15}). 

In 2011 Parkkonen and Paulin \cite{PaPa11} defined and studied quadratic Lagrange spectra 
by considering the approximation by elements of the orbit of a given real quadratic
irrational number for the action by homographies and anti-homographies of
$PSL_2(\Z)$ on $\R \cup \{\infty\}$. These spectra were further investigated 
by Bugeaud \cite{Bu14}, Pejkovi\'c \cite{Pej16}, and Lin \cite{Lin18}. Among other results, the existence 
of an Hall's ray for every quadratic Lagrange spectrum has been established in \cite{Lin18}. 
Subsequently, Parkkonen and Paulin \cite{PaPa18} defined and studied quadratic Lagrange 
spectra in completion of function fields over finite fields with respect to the ansolute values defined by discrete 
valuations. This setting includes the special case of the field of rational fractions and its 
valuation at infinity, which was given special attention in \cite{PaPa18} and is 
studied in the present paper. 

Let $q$ be a positive power of a prime number and 
$\kk = \FF_q$ denote the finite field of order $q$.  
Let $R=\FF_q[Y]$, $K=\FF_q(Y)$, and $\wh{K}=
\FF_q((Y^{-1}))$ be, respectively, the ring of polynomials in one
variable $Y$ over $\FF_q$, the field of rational functions in $Y$ over
$\FF_q$, and the field of formal power series in $Y^{-1}$ over
$\FF_q$.  Then, $\wh{K}$ is a nonarchimedean local field, the
completion of $K$ with respect to the
absolute value defined by $|{P \over Q}|=q^{\deg P-\deg Q}$ for all 
$P,Q$ in $R \setminus \{0\}$. Also, sometimes, it is convenient to use the associated valuation $\v$ 
defined for a nonzero element $f$ in $\wh{K}$ by
$$
\v (f) = - {\log |f| \over  \log q}.
$$ 

We stress that $\wh{K}$ is not algebraically closed (for instance, the polynomial $X^2 - Y$ has no 
roots in $\wh{K}$). Let
$$
\QI=\{f\in \wh{K}\;:\;[K(f):K]=2\}
$$ 
be the set of power series in $\wh{K}$ which are quadratic over $K$.  Given $f$ in 
$\wh{K} \setminus K$, it is well known that $f$ is in $\QI$ if and only if its
continued fraction expansion is eventually periodic. The projective action
of ${\rm PGL}_2(R)$ on $\wh K\cup \{\infty\}$ preserves
$\QI$, keeping the periodic part of the continued fraction expansions
unchanged, up to cyclic permutation and invertible elements; see Lemma 2.5 below. We refer
for instance to \cite{Las00,Schm00,Pau02} for
background on the above notions.
  
In odd characteristic, the quadratic polynomials which are irreducible over $R$ are separable. This is 
not the case in characteristic $2$, where there exist quadratic polynomials in $R[X]$ 
which are irreducible over $R$
and not separable. These polynomials are precisely the polynomials 
of the form $A X^2 + C$, with $A, C$ in $R$, whose (double) roots are either in $K$ or in 
$\FF_q((Y^{-1/2})) \setminus \wh{K}$. Consequently, every quadratic power series in $\wh{K}$ is a 
root of an irreducible polynomial of the form $A X^2 + B X + C$, where $A, B, C$
are in $R$ and $B$ is nonzero. 
This short discussion explains that, in characteristic $2$, the 
Galois conjugate $\alpha^\sigma$ of a quadratic power series $\alpha$ 
in $\wh{K}$ is well defined and that $\alpha - \alpha^\sigma = \alpha + \alpha^\sigma$ is nonzero.

Let $\alpha$ be in $\QI$ and $\alpha^\sigma$ in $\QI$
denote its Galois conjugate over $K$.  The {\it complexity}
$h(\alpha)= {1 \over |\alpha-\alpha^\sigma|}$ of $\alpha$ was
introduced in \cite{HePa10} and studied in Section 17.2
of \cite{BPP}. 
It measures the size of $\alpha$, in the same way as $\max\{|p|, |q|\}$ measures the 
size of the rational number $p/q$, where $p, q$ are nonzero coprime integers.
The complexity $h(\alpha)$ can be expressed in terms of the 
continued fraction expansion of $\alpha$; see Lemma 2.2. 
Let
$$
\Theta_\alpha={\rm PGL}_2(R) \cdot\{\alpha,\alpha^\sigma\}
$$
be the union of the orbits of $\alpha$ and $\alpha^\sigma$ under the
projective action of ${\rm PGL}_2(R)$. Given $f$ in $\wh K \setminus (K\cup
\Theta_\alpha)$, Parkkonen and Paulin \cite{PaPa18}
introduced the {\it quadratic approximation constant}
of $f$, defined by
$$
c_\alpha(f)=\liminf_{\beta\in\Theta_\alpha,\;|\beta-\beta^\sigma|\to 0}
\, {|f-\beta| \over |\beta-\beta^\sigma|}\;,
$$
and they studied the {\it quadratic Lagrange spectrum} of $\alpha$, defined by 
$$
\Sp(\alpha)=\{c_\alpha(f) : f \in\wh K \setminus (K\cup \Theta_\alpha)\}\;.
$$ 
Note that the quadratic Lagrange spectrum of $\alpha$ 
is contained in $q^\Z \cup \{0,+\infty\}$, thus, it is countable. It follows
from Theorem 1.6 of \cite{HePa10} that if $m_{\wh K}$ is a Haar measure
on the locally compact additive group of $\wh K$, then for $m_{\wh  K}$-almost every 
$f$ in $\wh K$, we have $c_\alpha(f)=0$. Hence in
particular, $0$ is in $\Sp(\alpha)$ and the quadratic Lagrange spectrum is
therefore closed. Parkkonen and Paulin \cite{PaPa18} proved that it is
bounded and defined the {\it (quadratic) Hurwitz
constant} of $\alpha$, denoted by $\Hw (\alpha)$, by
$$
\Hw (\alpha) := \max\Sp(\alpha)\;\in q^{\Z}.
$$
They obtained several results on $\Hw (\alpha)$ and on $\Sp(\alpha)$, including the 
existence of an Hall's ray, for every $\alpha$. 
Parkkonen and Paulin \cite{PaPa18} established nonarchimedean analogues of the results 
obtained in \cite{Bu14,Pej16,Lin18} by using geometric methods of group actions on Bruhat--Tits trees. 
In the present paper, 
we reprove many of their results by means of the theory of continued 
fractions in power series fields and, in addition, we establish several new results. 
We stress that, unlike in \cite{PaPa18}, we do not 
assume that the prime power $q$ is odd (however, Fr\'ed\'eric Paulin informed me that their 
approach also works well in characteristic $2$). 

In particular, we give alternative proofs of the following two theorems 
highlighted in~\cite{PaPa18}.

\proclaim Theorem 1.1 (Parkkonen and Paulin \cite{PaPa18}).
Let $\alpha$ be a quadratic power series in $\wh{K}$.
\smallskip \noindent 
(1) (Upper bound) Its quadratic Hurwitz constant satisfies
$\Hw(\alpha)\leq q^{-2}$.
\smallskip \noindent 
(2) (Hall's ray) There exists an integer $m_\alpha$ such that, 
for every integer $n$ with $n > m_\alpha$, the real number $q^{-n}$ belongs to $\Sp(\alpha)$.

\proclaim Theorem 1.2 (Parkkonen and Paulin \cite{PaPa18}). 
The Hurwitz constant of any
quadratic power series in $\wh{K}$, whose continued fraction
expansion is eventually periodic with a period of length at most $q-1$, is equal to
$q^{-2}$.

Examples of quadratic power series for which the quadratic
Lagrange spectrum coincides with its Hall ray are given in \cite{PaPa18}
and in Theorem 5.1. In particular, we have
$$
\Sp( [0; Y, Y, \ldots ] )=\{0\}\cup\{q^{-n-2} : n\in\Z_{\ge 0}\}.   \eqno (1.1) 
$$

Our approach shows that, to determine the quadratic spectrum of a quadratic power series $\alpha$, it is 
sufficient to compute the quadratic approximation constants of the quadratic power series
not in $\Theta_\alpha$. 

\proclaim Theorem 1.3. 
Let $\alpha$ be a quadratic power series in $\wh{K}$ and $q^{-m}$ a nonzero element of its spectrum. 
Then, there exists a quadratic power series $f$ such that $c_{\alpha} (f) = q^{-m}$. 

The proof of Theorem 1.3 shows that, if $d$ denotes the 
maximal degree of a partial quotient in the 
periodic part of the continued fraction expansion of $\alpha$, then one can in addition impose 
that all the partial quotients of $f$ have degree at most $d+1$. Moreover, one can also impose 
that the length of the period of the continued fraction expansion of $\alpha$ is at most equal 
to $2 + {m \over 2}$. Consequently, an integer $m$ being given, it is sufficient to compute 
$c_\alpha (f)$ for $f$ being in an explictly given finite set in order to determine whether $q^{-m}$ 
is or not in the spectrum of $\alpha$. Since a suitable value for the integer $m_\alpha$ in 
Theorem 1.1 (2) can be given explicitly, all this shows that a finite amount of 
computation is sufficient to determine exactly the set $\Sp (\alpha)$. 

\medskip

Proposition 4.9 of \cite{PaPa18} asserts that the function $\Hw$ takes arbitrarily small 
positive values. We establish that it can take every admissible value.

\proclaim Theorem 1.4. 
For every $m \ge 2$, there exists $\alpha$ in $\QI$ such that $\Hw(\alpha) = q^{-m}$.

Parkkonen and Paulin \cite{PaPa18} gave explicit examples of classes of quadratic
power series whose quadratic Lagrange spectrum does not coincide with
its Hall ray, in other words, which have at least one {\it gap} in their
spectrum. We go slightly further and establish that the number of gaps can be prescribed.

\proclaim Theorem 1.5. 
For any positive integer $\ell$, there exist quadratic power series in $\wh{K}$ 
whose Lagrange spectrum has exactly $\ell$ gaps.

Throughout, for $a_0$ in $R$ and $a_1, \ldots , a_{r+s}$ nonconstant 
polynomials in $R$, we use the notation
$$
[a_0; a_1, a_2 \ldots , a_r, \overline{a_{r+1}, \ldots , a_{r+s}}] 
:=  a_0 + {1 \over \disp a_1 +
{\strut 1 \over \disp a_2 + {\strut 1 \over \disp \ldots}}}
$$
to indicate that the block of partial quotients $a_{r+1}, \ldots , a_{r+s}$
is repeated infinitely many times.

We recall that an irrational power series $\alpha$ is quadratic if and only if 
its continued fraction expansion is ultimately periodic, that is, of the
form
$$
\alpha = [a_0; a_1, \ldots , a_r, \overline{b_{1}, \ldots , b_{s}}].   \eqno (1.2) 
$$
When we express $\alpha$ as in (1.2) we tacitly assume that $s$ is minimal and that 
$a_r \not= b_s$.  We call $b_1, \ldots , b_s$ the shortest periodic part
in the continued fraction expansion of $\alpha$.

It does not seem to be easy to find an expression of the Hurwitz constant of a
quadratic power series $\alpha$. 
Inspired by the results obtained in \cite{Bu14,Pej16,Lin18}, the authors of \cite{PaPa18} conjectured
that
$$
\Hw (\alpha) = \max \bigl\{\limsup_{P\in\FF_q[X],\;\deg P \to +\infty} 
c_\alpha([0; \overline{P}\;]),\;\;\max_{P\in\FF_q[X],\;\deg P=1} 
c_\alpha([0; \overline{P}\;])\bigr\},
$$ 
and established this formula in some special cases. 
At the end of Section 5, we give an example showing that this 
conjecture does not hold. Actually, we feel that there is no simple 
formula for $\Hw(\alpha)$.

The present paper is organized as follows. In Section 2, we gather several results on continued fractions 
in power series fields and apply them in Section 3 to prove Theorem 1.3. 
Sections 4 and 5 are devoted to the proofs of our further results.

\vskip 6mm

\centerline{\bf 2. Auxiliary lemmas on continued fractions in power series fields}

\vskip 6mm

We assume that the reader is familiar with the classical theory of continued fractions of real numbers.
Good references include \cite{Per,HW} and \cite{Las00,Schm00} for the case of power series.

Our first lemma is an analogue for quadratic power series 
of a theorem of Galois.

\proclaim Lemma 2.1.
Let $s \ge 1$ be an integer and $b_1, \ldots , b_s$ nonconstant polynomials in $R$. 
The Galois conjugate of the quadratic power series
$$
\tau := [b_{1}; \overline{b_{2}, \ldots , b_{s}, b_{1}}]
$$
is the power series
$$
\tau^{\sigma} = - [0; \overline{b_{s}, \ldots , b_{2}, b_{1}}].
$$

\pro
Define
$$
{p_s \over q_s} := [0; b_1, b_2, \ldots b_s], \quad 
{p_{s-1} \over q_{s-1}} := [0; b_1, b_2, \ldots b_{s-1}]
$$
and
$$
{p'_s \over q'_s} := [0; b_s, b_{s-1}, \ldots b_1], \quad 
{p'_{s-1} \over q'_{s-1}} := [0; b_s, b_{s-1}, \ldots b_2].
$$
Then, $\tau$ satisfies
$$
{1 \over \tau} = {p_s \tau + p_{s-1} \over q_s \tau + q_{s-1}},
$$
hence,
$$
p_s \tau^2 + (p_{s-1} - q_s) \tau - q_{s-1} = 0.    \eqno (2.1)
$$
Likewise, $\tau' :=  [0; \overline{b_{s}, \ldots , b_{2}, b_{1}}]$ satisfies
$$
\tau' = {(p'_s / \tau') + p'_{s-1} \over (q'_s / \tau') + q'_{s-1}} 
= {p'_{s-1} \tau' + p'_s \over q'_{s-1} \tau' + q'_s},
$$
hence,
$$
q'_{s-1} (\tau')^2 + (q'_s - p'_{s-1}) \tau' - p'_s = 0.  \eqno (2.2)
$$
The mirror formula (see e.g. \cite{Per}, page 32) gives us that 
$$
p'_s = q_{s-1}, \quad q'_s = q_s, \quad 
p'_{s-1} = p_{s-1}, \quad \hbox{and} \quad q'_{s-1} = p_s.
$$
Combined with (2.2), we obtain
$$
p_s (\tau')^2 + (q_s - p_{s-1}) \tau' - q_{s-1} = 0.  \eqno (2.3)
$$
Equalities (2.1) and (2.3) show 
that $\tau$ and $-\tau'$ are roots of the same quadratic polynomial. Since they are distinct, they are 
Galois conjugate. 
\cqfd

Our second lemma establishes that 
the quantity $h(\alpha) = |\alpha - \alpha^{\sigma}|^{-1}$ can be expressed in a 
simple way in terms of the continued fraction expansion
of the quadratic power series~$\alpha$.

\proclaim Lemma 2.2.
Let $\alpha$ be a quadratic power series with ultimately periodic
continued fraction expansion given by 
$$
\alpha = [a_0; a_1, \ldots , a_r, \overline{b_{1}, \ldots , b_{s}}], 
$$
where $s \ge 1$. Denote by $\alpha^{\sigma}$ its Galois conjugate.
If $a_r \not= b_{s}$, then we have
$$
h(\alpha) = |\alpha - \alpha^{\sigma}|^{-1} 
= q^{2 ( \sum_{i=1}^r \deg a_i ) - \deg a_r - \deg b_s +  \deg(a_r - b_s)}.
$$

Lemma 2.2 is the analogue of Lemma 2.1 of \cite{Bu14} (see also Lemma 6.1 of \cite{Bu12}).

\pro
By Lemma 2.1, the Galois conjugate of 
$$
\tau := [b_{1}; \overline{b_{2}, \ldots , b_{s}, b_{1}}]
$$
is the quadratic number
$$
\tau^{\sigma} = - [0; \overline{b_{s}, \ldots , b_{2}, b_{1}}].
$$
Let $(p_\ell/q_\ell)_{\ell \ge 1}$ denote the sequence of
convergents to $\alpha$. Since we have
$$
\alpha = {p_r \tau + p_{r-1} \over q_r  \tau + q_{r-1}}
\quad \hbox{and} \quad
\alpha^{\sigma}  = {p_r \tau^{\sigma} + p_{r-1} \over q_r \tau^{\sigma} + q_{r-1}},
$$
we get
$$
\eqalign{
|\alpha - \alpha^{\sigma} | & =    { |\tau - \tau^{\sigma}|  \over 
|q_r  \tau + q_{r-1}| \cdot |q_r \tau^{\sigma} + q_{r-1}| } \cr
& = q^{-2 \deg q_r} \, |\tau - \tau^{\sigma} | \,  \Bigl|\tau + {q_{r-1} \over q_r} \Bigr|^{-1}
\, \Big|\tau^{\sigma} + {q_{r-1} \over q_r}\Bigr|^{-1}. \cr} \eqno (2.4)
$$
Observe that
$$
|\tau - \tau^{\sigma}| = \Bigl|\tau + {q_{r-1} \over q_r} \Bigr| 
= q^{\deg b_1}    \eqno (2.5) 
$$
and
$$
\eqalign{
\Big|\tau^{\sigma} + {q_{r-1} \over q_r}\Bigr| & = 
\bigl| [0; a_r, a_{r-1}, \ldots ] - [ 0; \overline{b_{s}, \ldots , b_{2}, b_{1}}] \bigr| \cr
& = \Biggl| {[a_r; a_{r-1}, \ldots ] - [b_s ; \overline{b_{s-1}, \ldots , b_{1}, b_{s}}] \over 
[a_r; a_{r-1}, \ldots ]  \cdot [b_s ; \overline{b_{s-1}, \ldots , b_{1}, b_{s}}] } \Biggr|  \cr
& = q^{\deg(a_r - b_s) - \deg a_r - \deg b_s}. \cr}
$$
By (2.4) and (2.5), this completes the proof of the lemma, 
since $\deg q_r =  \sum_{i=1}^r \deg a_i$. 
\cqfd

Our third auxiliary lemma is the analogue of Lemma 2.2 from \cite{Bu14}.

\proclaim Lemma 2.3. 
Let $\alpha = [0; a_1, a_2, \ldots]$ and $\beta =
[0; b_1, b_2, \ldots]$ be power series in $\wh{K}$. 
Assume that there exists a nonnegative integer $n$ such that  
$a_i = b_i$ for any $i=1, \ldots, n$ and $a_{n+1}\not=b_{n+1}$. 
Then, we have
$$
|\alpha - \beta| = q^{-2 (  \sum_{i=1}^n \deg a_i ) - \deg a_{n+1} - \deg b_{n+1} 
+ \deg (a_{n+1} - b_{n+1}) }. 
$$

\pro 
Set $\alpha'=[a_{n+1}; a_{n+2},\ldots]$ and 
$\beta'=[b_{n+1}; b_{n+2},\ldots]$. 
Let $(p_\ell / q_\ell)_{\ell \ge 1}$ denote 
the sequence of convergents to $\beta$. Since $a_{n+1}\not= b_{n+1}$ 
and the first $n$ partial quotients of $\alpha$ and
$\beta$ are assumed to be the same, we get 
$$
\alpha={p_n\alpha'+p_{n-1}\over q_n\alpha'+q_{n-1}} \quad {\rm and} \quad
\beta={p_n\beta'+p_{n-1}\over q_n\beta'+q_{n-1}},
$$
thus,
$$
\vert\alpha-\beta\vert=
\left\vert{p_n\alpha'+p_{n-1}\over q_n\alpha'+q_{n-1}}-
{p_n\beta'+p_{n-1}\over q_n\beta'+q_{n-1}}\right\vert
= \left\vert{ \alpha'-\beta'
\over (q_n\alpha'+q_{n-1})(q_n\beta'+q_{n-1})}\right\vert\cdot   
$$
Since $\deg q_n = \deg a_1 + \cdots + \deg a_n$, 
$$
| q_n\alpha'+q_{n-1} | = q^{-\deg q_n - \deg a_{n+1}}, \quad 
| q_n\beta'+q_{n-1} | = q^{-\deg q_n - \deg b_{n+1}},
$$
and
$$
| \alpha'-\beta' | = q^{\deg (a_{n+1} - b_{n+1})}, 
$$
this proves the theorem. 
\cqfd

We display an easy consequence of Lemma 2.3. 
Below and in the next sections, 
it is convenient to take a point of view from combinatorics on words. 
For an integer $k \ge 1$, let ${\cal A}_{\le k}$ (resp., ${\cal A}_{= k}$) 
denote the set of all nonconstant polynomials 
in $R$ of degree at most equal to $k$ (resp., equal to $k$). Set 
$$
{\cal A} := \bigcup_{k \ge 1} \, {\cal A}_{\le k} =\bigcup_{k \ge 1} \, {\cal A}_{= k}.
$$ 
If $a_1 \ldots a_r$ is a finite word over ${\cal A}$, then $(a_1 \ldots a_r)^{\infty}$ denotes the 
infinite word obtained by concatenating infinitely many copies of $a_1 \ldots a_r$.

\proclaim Corollary 2.4.
Let
$$
\tau = [0; \overline{b_{1}, b_2 \ldots , b_{s}}]
$$
be a quadratic power series whose shortest periodic part is $b_1, \ldots , b_s$. Let
$$
f = [a_0; a_1, a_2, \ldots]
$$
be an irrational real number not in $\Theta_{\tau}$.
For a positive integer $r$ such that $a_r \not= b_s$, set
$$
\alpha_r := [a_0; a_1, a_2, \ldots , a_r, \overline{b_{1}, \ldots , b_{s-1}, b_{s}}].
$$ 
If $a_{r+1} \not = b_1$, then 
$$
\eqalign{
\v \biggl( { |f - \alpha_r|  \over |\alpha_r - \alpha_r^{\sigma} |} \biggr)  
=   \deg a_r  & + \deg b_s - \deg (a_r - b_s)  \cr
& + \deg a_{r+1} + \deg b_{1} - \deg (a_{r+1} - b_{1}).  \cr}   \eqno (2.6)
$$
If $a_{r+1} = b_1$, then 
let $t$ be the largest integer such that the word $a_{r+1} \ldots a_{r+t}$ coincide with the 
prefix of length $t$ of the infinite word $(b_{1} \ldots b_{s})^{\infty}$. 
Let $s_0$ be defined by $a_{r+t} = b_{s_0}$ and put $s' = s_0 + 1$ if $s_0 < s$ and 
$s' = 1$ otherwise. Then, we have 
$$
\eqalign{ 
\v \biggl( { |f - \alpha_r|  \over |\alpha_r - \alpha_r^{\sigma} |} \biggr)  
=  2 \bigl( \, \sum_{j=1}^{t} \deg a_{r+j} \, \bigr) & + \deg a_r + \deg b_s - \deg (a_r - b_s)  \cr
& + \deg a_{r+t+1} + \deg b_{s'} - \deg (a_{r+t+1} - b_{s'})  .   \cr}  
$$
In particular, in both cases, we have 
$$
{ |f - \alpha_r|  \over |\alpha_r - \alpha_r^{\sigma} |} \le q^{-2}.  \eqno (2.7) 
$$
Furthermore, putting $t=0$ and $s'=1$ if $a_{r+1} \not = b_1$, we obtain that, if 
$\deg a_r \not= \deg b_s$ and $\deg a_{r+t+1} \not= \deg b_{s'}$, then 
$$
\eqalign{
\v \biggl( { |f - \alpha_r|  \over |\alpha_r - \alpha_r^{\sigma} |} \biggr) 
=   2 \bigl( \, \sum_{j=1}^{t} \deg a_{r+j} \,  \bigr) & + \min \{\deg a_r, \deg b_s\}  \cr
& + \min \{\deg a_{r+t+1}, \deg b_{s'} \}. \cr}     \eqno (2.8)
$$

\pro
This follows directly from Lemmas 2.2 and 2.3. 
\cqfd 

It remains for us to describe the orbit of a quadratic power series under the action 
of ${\rm PGL}_2(R)$. 
For the real analogue, that is, for characterising the orbit of an irrational number under the action 
of ${\rm SL}_2 (\Z)$, we used in \cite{Bu14} a classical theorem of Serret (see \cite{Per}, page 65), 
which asserts that the tails of the continued fraction expansions of 
two irrational real numbers $\alpha, \beta$ coincide
if and only if there exist integers $a, b, c, d$ with $ad-bc = \pm 1$
such that
$$
\alpha = {a \beta + b \over c \beta + d}.
$$

In the present context, we make use of the version of Serret's theorem established by 
Schmidt (Theorem 1 of \cite{Schm00}). 
Before stating it, let us observe that, for any irrational power series 
$\alpha := [a_0; a_1, a_2, \ldots]$ in $\wh K$ and any $a$ in $\kk^{\times}$, we have
$$
a \alpha = [a a_0; a^{-1} a_1, a a_2, a^{-1} a_3, \ldots ].
$$
Furthermore, two power series $\alpha, \beta$ in $\wh K$ are called {\it equivalent} if there exist 
$a, b, c, d$ in $R$ with $a d - b c$ in $\kk^{\times}$ and such that 
$$
\alpha = {a \beta + b \over c \beta + d}.
$$

\proclaim Lemma 2.5. 
Two irrational power series $\alpha := [a_0; a_1, a_2, \ldots]$ and 
$\beta := [b_0; b_1, b_2, \ldots ]$ in $\wh K$ are equivalent if and only if there exist nonnegative 
integers $m, n$ and an element $a$ in $\kk^{\times}$ such that
$$
\alpha := [a_0; a_1, a_2, \ldots , a_{n-1}, a_n, a_{n+1}, a_{n+2}, \ldots]
$$
and
$$
\beta := [b_0; b_1, b_2, \ldots , b_{m-1}, a a_n, a^{-1} a_{n+1}, a a_{n+2}, \ldots].
$$

\pro 
This is Theorem 1 of \cite{Schm00}. It is also proved in 
Section IV.3 of \cite{Math70}.
\cqfd

\vskip 6mm

\goodbreak 

\centerline{\bf 3. An equivalent formulation for $\Sp (\alpha)$}

\vskip 6mm

Throughout this section, we fix a quadratic power series $\alpha$. 
Let $b_1, \ldots , b_s$ be the (shortest) periodic part in its
continued fraction expansion and set
$$
\tau = [b_{1}; \overline{b_{2}, \ldots , b_{s}, b_{1}}].
$$
For $j = 1, \ldots , s$ and $a$ in $\kk^{\times}$, set
$$
\tau_{j, a} := [a b_{j} ;  \overline{a^{-1} b_{j+1}, \ldots , a^{(-1)^{j-1}} b_{j-1}, a^{(-1)^{j}} b_{j}, 
a^{(-1)^{j+1}} b_{j+1},  \ldots , a^{-1} b_{j-1}, a b_{j}}]
$$
and
$$
\tau'_{j, a}  = [a^{-1} b_{j-1} ; \overline{a b_{j - 2}, a^{-1} b_{j-3}, \ldots , a^{(-1)^{j}} b_{j}, 
a^{(-1)^{j-1}} b_{j-1}, a^{(-1)^{j-2}}  b_{j-2}, 
\ldots , a b_{j}, a^{-1} b_{j-1}}]. 
$$
Here and below, the indices are always understood modulo $s$. 
Observe that $\tau = \tau_{1, 1}$ and, for $a$ in $\kk^{\times}$, we have 
$$
\Theta_{\tau} = \Theta_{\tau_{1, a}}  = \ldots = \Theta_{\tau_{s, a}} = \Theta_{\tau'_{1,a}} 
= \ldots = \Theta_{\tau'_{s, a}} = \Theta_{\alpha}. 
$$
Observe also that the length of the shortest periodic part of $\tau_{j, a}$ divides $2s$. 
Furthermore, by Lemma 2.1, we have
$$
\tau_{j, a}^{\sigma} = - [0; \overline{a^{-1} b_{j - 1}, a b_{j-2}, \ldots , a^{-1} b_{j+1}, a b_{j}}] = - 1/\tau'_{j, a},
$$
for $j=1, \ldots , s$.

Let
$$
f = [0; a_1, a_2, \ldots]
$$
be an irrational power series not in $\Theta_{\tau}$,
which we wish to approximate by power series from $\Theta_{\tau}$. 
A trivial way to do this consists in keeping the $r$ first partial quotients 
of $f$ and putting then the sequence of partial quotients of one 
of the power series $\tau_{j, a}$ or $\tau'_{j, a}$, with $1 \le j \le s$ and $a$ in $\kk^\times$ 
(recall that $\Theta_{\tau}$ 
is equal to ${\rm PGL}_2(R) \cdot\{\alpha,\alpha^\sigma\}$). For instance, for $r \ge 1$
and $j = 1, \ldots , s$, the quadratic power series
$$
\alpha_{r, j, a} := [0; a_1, a_2, \ldots , a_r, \overline{a b_{j}, a^{-1} b_{j+1}, \ldots , a^{-1} b_{j-1}}]
= [0; a_1, a_2, \ldots , a_r, \tau_{j, a}]
$$
and
$$
\alpha'_{r, j, a} := [0; a_1, a_2, \ldots , a_r, \overline{a^{-1} b_{j-1}, a b_{j-2}, \ldots , a b_{j}}]
= [0; a_1, a_2, \ldots , a_r, \tau'_{j, a}]
$$
are quite good approximations to $f$ in $\Theta_{\tau}$ and
$$
\ell_{j, a} (f) := \liminf_{r \to + \infty} \, |f - \alpha_{r,j, a} | \cdot h(\alpha_{r,j, a}),  \quad
\ell'_{j, a} (f) := \liminf_{r \to + \infty} \, |f - \alpha'_{r,j, a} | \cdot h(\alpha'_{r,j, a })   
$$
are greater than or equal to $c_{\tau} (f)$, thus
$$
c_{\tau} (f) \le \min_{1 \le j \le s} \, \min_{a \in \kk^{\times}} \, 
\min\{ \ell_{j, a} (f), \ell'_{j, a} (f)\}.
\eqno (3.1)
$$
Unlike what happens in the real case, we do have equality in (3.1). 


\proclaim Lemma 3.1. 
Under the above notation, we have
$$
c_{\tau} (f) = \min_{1 \le j \le s} \, \min_{a \in \kk^{\times}} \,  \min\{ \ell_{j, a} (f), \ell'_{j, a} (f)\}.
$$

The analogue of Lemma 3.1 does not hold in the real case, see e.g. Section 3.6 of \cite{Bu14}. 
Lemma 3.1 shows that the power series case is simpler than its real analogue. 

\pro
We have to estimate the
quantities $|f - \alpha| \cdot h(\alpha)$ for quadratic power series $\alpha$ of the form
$$
\zeta_{j, a} := [0; a_1, a_2, \ldots , a_r, c_1, \ldots , c_t, \tau_{j, a}]
\quad \hbox{and} \quad
\zeta'_{j, a} := [0; a_1, a_2, \ldots , a_r, c'_1, \ldots , c'_t, \tau'_{j, a}],
$$
where $1 \le j \le s$, $t \ge 1$, $a$ in $\kk^\times$, and $c_1, \ldots , c_t, c'_1, \ldots , c'_t$ are nonconstant 
polynomials in $R$, with $c_1 \not= a_{r+1}$, $c_t \not= b_{j-1}$, 
$c'_1 \not= a_{r+1}$, and $c'_t \not= b_j$. For simplicity, we only treat 
the case of $\tau_{j, 1}$, the other cases being completely analogous. 

It follows from Lemmas 2.2 and 2.3 that 
$$
\eqalign{
\v \biggl( {|f - \zeta_{j, 1}| \over |\zeta_{j, 1} - \zeta_{j, 1}^{(\sigma)}|}  \biggr) 
& =   \deg a_{r+1} + \deg c_1 - \deg (a_{r+1} - c_1) \cr 
& \hskip 7mm - 2 \sum_{j=1}^t \deg c_j  +
\deg c_t + \deg b_{j-1} - \deg (c_t - b_{j-1}),  \cr} 
\eqno (3.2) 
$$
where $j-1$ is understood modulo $s$. 
If $t \ge 2$, then the right hand side of (3.2) is
$$
\le - \deg (a_{r+1} - c_1) + \deg a_{r+1} - \deg c_1 - \deg c_t - \deg (c_t - b_{j-1}) + \deg b_{j-1} \le 0. 
$$
If $t=1$, then the right hand side of (3.2) is equal to 
$$
- \deg (a_{r+1} - c_1) + \deg a_{r+1}  - \deg (c_1 - b_{j-1}) + \deg b_{j-1}. 
$$
Since $\deg (P_1 + P_2) \le \max\{\deg P_1, \deg P_2\}$ holds for all polynomials $P_1, P_2$, 
the latter quantity is 
$$
 \le - \deg (a_{r+1} - b_{j-1}) + \deg a_{r+1} + \deg b_{j-1}. 
$$
Recalling that 
$$
\alpha_{r, j-1, 1} := [0; a_1, a_2, \ldots , a_r, \tau_{j-1, 1}], 
$$
it follows from (2.6) that 
$$
{|f - \alpha_{r, j-1, 1}| \over |\alpha_{r, j-1, 1} - \alpha_{r, j-1, 1}^{(\sigma)}|} \le q^{- \deg a_{r+1} - \deg b_{j-1}
+  \deg (a_{r+1} - b_{j-1}) }
\le {|f - \zeta_{j, 1}| \over |\zeta_{j, 1} - \zeta_{j, 1}^{(\sigma)}|}.
$$
Consequently, it is sufficient to restrict our attention to the approximants of the form 
$\alpha_{r,j,a}$ and $\alpha'_{r,j,a}$ in order to compute $c_{\tau} (f)$. This proves the lemma. 
\cqfd

\noindent {\bf Notation. } 
Let $W = w_1 \ldots w_h$ with $h \ge 1$ denote a finite word over the alphabet ${\cal A}$. 
Then, $f_W$ denotes the quadratic power series with purely periodic continued fraction expansion
of period $W$, that is, 
$$
f_W = [0; \overline{w_1, \ldots , w_h} ].
$$

\medskip

An important consequence of Lemma 3.1 is the fact 
that the spectrum of $\alpha$ is determined by the set of values 
taken by the function $c_\alpha$ at quadratic power series. This is 
precisely the content of Theorem 1.3.


\medskip

\noi {\it Proof of Theorem 1.3. } 
Let $f$ be in $\wh K$ and assume that there exists a positive integer $m$ such that 
$c_\alpha (f) = q^{-m}$. Then, there exists $j$ in $\{1, \ldots , s\}$ and $a$ in $\kk^{\times}$ such that 
$\ell_{j, a} (f) = q^{-m}$ or $\ell'_{j, a} (f) = q^{-m}$. 
Without any loss of generality, we may assume that $a = 1$ 
and $\ell_{j, 1} (f) = q^{-m}$. Consequently, there 
are arbitrarily large integers $r$ such that 
$$
\alpha_{r,j, 1} := [0 ; a_1, a_2, \ldots , a_r, \tau_{j, 1}] 
$$
satisfies
$$
|f - \alpha_{r,j, 1} | \cdot h(\alpha_{r,j, 1}) = q^{-m}.
$$
There are polynomials $P_1, P_2$ and a word $W_0$
(possibly empty) which is a factor of $(b_1 b_2 \ldots b_s)^{\infty}$ 
such that there exist $t \ge 1$ and arbitrarily large integers $r$ with
$$
a_r = P_1, \quad a_{r+1} \ldots a_{r+t} = W_0, \quad a_{r+t+1} = P_2,
$$
and
$$
|f - \alpha_{r,j, 1} | \cdot h(\alpha_{r,j, 1}) = q^{-m}.
$$
Setting $W = P_1 W_0 P_2$, we see that $c_\alpha (f_W) \le q^{-m}$. The inequality cannot be strict,
since otherwise we would have $c_\alpha (f) < q^{-m}$. Consequently, $c_\alpha (f_W)$ is equal 
to $q^{-m}$ and the theorem is proved. 
\cqfd

\vskip 6mm

\goodbreak 

\centerline{\bf 4. First results on $\Sp(\alpha)$ for an arbitrary $\alpha$}

\vskip 6mm

We display several immediate consequences of the preceding lemmas. 
Our first result is a reformulation of the first assertion of Theorem 1.1. 

\proclaim Theorem 4.1.
For every quadratic power series
$\alpha$ in $\wh K$, the spectrum $\Sp(\alpha)$ satisfies
$$
\Sp(\alpha) \subset \{0\}\cup\{q^{-n-2} : n \in \Z_{\ge 0} \}.
$$  

\pro
This follows from (2.7) and Lemma 3.1.
\cqfd

Following \cite{PaPa18}, for every power series $f$ in $\wh K \setminus K$, set 
$$
\eqalign{
M(f)  & :=\limsup_{k\to+\infty}\;\deg a_k\;\geq 1,\cr
M_2(f) & :=\limsup_{k\to+\infty}\big(\deg a_k+\deg a_{k+1}\big)\;\geq 2,\cr
m(f) & :=\liminf_{k\to+\infty}\;\deg a_k\;\geq 1. \cr}
$$
Corollary 2.4 allows us to reprove Lemma 4.4 and Corollaries 4.6 and 4.7 of \cite{PaPa18}.

Throughout the end of this section, we keep the notation of Section 3. 
We denote by $b_1, \ldots , b_s$ the (shortest) periodic part in the 
continued fraction expansion of a power series 
$\alpha$ in $\QI$ and we define $\tau$ and $\tau_{j, a}$
for $j = 1, \ldots , s$ and $a$ in $\kk^\times$ as in Section 3. 

\proclaim Proposition 4.2. 
Let $\alpha$ be in $\QI$ and $f$ in $\wh K \setminus (K\cup\Theta_\alpha)$.
\smallskip \noindent 
(1) If $m(f)>M(\alpha)$, then $c_\alpha(f) = q^{-M_2(\alpha)}$. Consequently,  
$\Hw (\alpha) \ge q^{-M_2(\alpha)}$.
\smallskip\noindent 
(2) If $M(f)<m(\alpha)$, then $c_\alpha(f) =  q^{-M_2(f)}$. Consequently, the quadratic spectrum of $\alpha$ 
includes $q^{-2}, q^{-3}, \ldots , q^{- 2 m(\alpha) + 2}$. 
\smallskip\noindent 
(3) If $M(\alpha) = 1$ or $m (\alpha) \ge 2$, then $\Hw(\alpha) = q^{-2}$. 

\pro
Replacing if necessary $\alpha$ by $\tau_{j, 1}$ for a suitable $j$ in $\{1, \ldots , s\}$, 
we may assume that $\deg b_1 + \deg b_s = M_2 (\alpha)$.  
If $m(f)>M(\alpha)$, then, for $r$ large enough, the right hand side
of (2.8) is equal to $M_2 (\alpha)$. This proves the 
first assertion. 
For the second assertion, we consider the infinite sequence of integers $r$ 
such that $\deg a_r + \deg a_{r+1} = M_2(f)$ and use the 
quadratic power series $[a_0 ; a_1, \ldots , a_r ,  \overline{b_{1}, \ldots , b_{s-1}, b_{s}}]$
to approximate $f$.  Then, choosing $f= [0; \overline{Y^d, Y^{d'}}]$ for integers $d, d'$ with 
$1 \le d, d' < m(\alpha)$, we get the second part of (2).
The assertion (3) is an immediate consequence of the first two. 
\cqfd

We now confirm the existence of an Hall's ray in the quadratic Lagrange spectrum of an 
arbitrary quadratic power series. This establishes Theorem 1.1 (2). 

\proclaim Theorem 4.3. 
For every $\alpha$ in $\QI$,  there exists $m_{\alpha}$ such that $q^{-m}$ is in $\Sp (\alpha)$ 
for every integer $m \ge m_{\alpha}$. Denoting by $s$ the length of the periodic part of the continued 
fraction expansion of $\alpha$ and by $d$ the maximum
of the degrees of its partial quotients, an admissible value for $m_\alpha$ is $2 d (s + 1)$. 

\pro
Replacing if necessary $\alpha$ by $\tau_{j, 1}$ for a suitable $j$ in $\{1, \ldots , s\}$, we 
may assume that $d := \deg b_s$ is equal to 
the maximum of the degrees of $b_1, \ldots , b_s$. 
Define $b_m$ for $m > s$ by setting $b_m = b_j$, where $1 \le j \le s$ and $j$ 
is congruent to $m$ modulo $s$. 
Let $u$ be an integer with $u \ge s$. 
Let $P_1, P_2$ 
be nonconstant polynomials in $R$ such that $P_1 \not= b_s$ and $P_2 \not= b_{u+1}$. 
It follows from Corollary 2.4 that
$$
\eqalign{
-  & {\log c_\tau ( [0; \overline{ Y^{d+1}, Y^{d+1}, P_1, b_1, \ldots , b_u, P_2 } ]  ) \over \log q}  \cr
&    \ge 2  \sum_{i=1}^u \deg b_i  + \deg b_s + \deg P_1 - \deg (b_s - P_1) 
+ \deg b_{u+1} + \deg P_2 - \deg (b_{u+1} - P_2),  \cr}  
$$
with equality since $u \ge s$ and our assumption that $P_1 \not= b_s$ and $P_2 \not= b_{u+1}$
guarantees that neither $P_1 b_1 \ldots b_s$, nor 
$b_2 \ldots b_s P_2$ is a factor of the period of some $\tau_{j, a}$. 
Some condition on $u$ is indeed necessary: it may happen 
that $P_1$ is one of $b_1, \ldots , b_{s-1}$, say $P_1 = b_\ell$, and the word $b_{\ell} b_1 b_2$ 
is a factor of $b_1 \ldots b_s$. For $u \ge s$, 

Selecting $P_1 = P_2 = Y^{d+1}$ and recalling that $\deg b_s = d$, this shows at once that 
$$
q^{- 2  \sum_{i=1}^u \deg b_i  - d - \deg b_{u+1}} 
$$
is in $\Sp (\alpha)$.

To establish the theorem, it is sufficient to show that, with suitable choices of $P_1$ and $P_2$, the quantity 
$$
\deg b_s + \deg P_1 - \deg (b_s - P_1) + \deg b_{u+1} + \deg P_2 - \deg (b_{u+1} - P_2)
$$
takes every integer value between $d +\deg b_{u+1}$ and 
$d + 2 \deg b_{u+1} + \deg b_{u+2}$. We proceed as follows. 
For $k = 0, \ldots , d$, there exists a polynomial $P_{1, k}$ of degree $d$ such that 
$\deg (P_{1, k} - b_s) = d - k$. 
For $h = 0, \ldots , \deg b_{u+1}$, there exists a polynomial $P_{2, h}$ 
of degree $\deg b_{u+1}$ such that 
$\deg (P_{2, h} - b_{u+1}) = \deg b_{u+1} - h$. Then, we have 
$$
\eqalign{
&\hskip 5mm  \deg b_s + \deg P_{1,k} - \deg (b_s - P_{1,k}) 
+ \deg b_{u+1} + \deg P_{2, h} - \deg (b_{u+1} - P_{2,h})  \cr
& = d + \deg b_{u+1} + k + h, \cr} 
$$
which takes all values between $d + \deg b_{u+1}$ and $d + 2 \deg b_{u+1} + \deg b_{u+2}$, 
since $d = \deg b_s \ge \deg b_{u+2}$. 
This shows that every rational number of the form $q^{-m}$, with $m$ 
an integer at least equal to $2  \sum_{i=1}^s \deg b_i  + d + \deg b_{s+1}$, is in $\Sp(\alpha)$. 
Consequently, a suitable value for $m_{\alpha}$ 
is given by $2  \sum_{i=1}^s \deg b_i  + d + \deg b_{u+1}$, which is at most equal to $2 d (s + 1)$. 
\cqfd

\vskip 6mm

\goodbreak 

\centerline{\bf 5. Gaps in the spectra and further results}

\vskip 6mm

We begin with an alternative proof of Theorem 4.11 of \cite{PaPa18} and establish (1.1). 

\proclaim Theorem 5.1 (Parkkonen and Paulin \cite{PaPa18}).  
For every polynomial $P$ in $R$ of degree $1$, we have 
$$
\Sp([0; \overline{P}\,])=\{0\}\cup\{q^{-n}:n\in\Z_{\ge 2}\}. 
$$

\pro
Set $\alpha := [0; \overline{P}\,]$. 
Let $m$ be a non-negative integer and set  
$$
g_m := [0 ; \overline{Y^2, Y^2, P, \ldots , P}] \quad \hbox{and} \quad 
h_m := [0 ; \overline{Y^2, Y^2, P+1, P, \ldots , P}],
$$ 
where $P$ is repeated $m$ times. Then, we check that $c_\alpha (g_m) = q^{-2m-2}$ 
and $c_\alpha (h_m) = q^{-2m-3}$. This shows that every $q^{-n}$ with $n \ge 2$ is in the spectrum 
of $\alpha$. Since we already observed that $0$ is in the spectrum of $\alpha$, this proves 
the theorem. 
\cqfd


We continue with an alternative proof of Proposition 4.8 of \cite{PaPa18}.

\proclaim Proposition 5.2 (Parkkonen and Paulin \cite{PaPa18}).
If $\alpha$ is in $\QI$ and the period of its continued fraction expansion 
contains no more than $q-2$ partial quotients of degree $1$, then 
$\Hw (\alpha) = q^{-2}$.

\pro
The argument is the same as in the proof of Proposition 4.8 of \cite{PaPa18}.
There exists a polynomial $P$ in $R$ of degree $1$ such that, for every 
partial quotient $b$ of degree $1$ of the period of $\alpha$, the polynomial 
$P - b$ is nonconstant. It then follows from Corollary 2.4 that 
$c_\alpha (f_P) = q^{-2}$. 
\cqfd 

\medskip

\noi {\it Proof of Theorem 1.2. } 
Theorem 1.2 follows from Proposition 5.2 if $k \le q-2$ or if $k = q-1$ and at least 
one partial quotient is of degree at least $2$. If $k = q-1$ and all the partial quotients 
are of degree $1$, then $M(\alpha) = 1$, $M_2 (\alpha) = 2$, and 
Theorem 1.2 follows from Proposition 4.2. 
\cqfd

\bigskip

\noindent {\it Proof of Theorem 1.4.} 
Let $k \ge 2$ be an integer. 
Consider a finite word $W$ over the alphabet ${\cal A}_{\le k}$ constructed by concatenating 
a copy of each different block of length $k-1$ over ${\cal A}_{\le k}$. 
The order is irrelevant.
Let $\alpha$ be in $\QI$ whose period is given by $W$. 

Let $f := [0; a_1, a_2, \ldots ]$ be in $\wh K \setminus (K\cup \Theta_\alpha)$. 
Let $h$ be the largest non-negative integer for which there are 
arbitrarily large integers $n$ such that the $h$ polynomials 
$a_n, a_{n+1}, \ldots , a_{n+h-1}$ are of degree at most $k$.
If $h \ge k-1$, then $c_{\alpha} (f) \le q^{-2 (k-1) - 2} = q^{-2k}$. 
If $h = k-2$ and $k \ge 3$, then there exists a polynomial $b$ of degree $k$ and infinitely many 
integers $n$ such that 
$a_n a_{n+1} \ldots  a_{n+h-1} b$ is a factor of $W$ and $\deg a_{n+h} > k$. This implies that 
$c_{\alpha} (f) \le q^{-2 (k-2) - 1 - k} = q^{- 3 k + 3}$. 
If $h < k-2$ or if $k=2$ and $h=0$, then there exists a polynomial $b$ of degree $k$ and infinitely many 
integers $n$ such that 
$b a_n a_{n+1} \ldots  a_{n+h-1} b$ is a factor of $W$, 
$\deg a_{n-1} > k$, and $\deg a_{n+h} > k$. This implies that 
$c_{\alpha} (f) \le q^{-2h - 2k} \le q^{- 2k}$. 
Noticing that 
$$
c_\alpha ( [0; \overline{ Y^{k+1}} ] ) = q^{-2k}, 
$$
all this shows that $\Hw(\alpha) = q^{-2k}$ for $k \ge 2$.

It remains for us to treat the case of $q^{-m}$ with $m$ odd. 
Consider a finite word $W'$ over the alphabet ${\cal A}_{\le k}$ constructed by concatenating 
a copy of each different block of length $k$ over ${\cal A}_{\le k}$. 
The order is irrelevant; however, for technical reasons, we assume that the last letter of $W'$ is $Y^k$. 
Let $\beta$ be in $\QI$ whose period is given by the word $W' Y^{k+1}$. 

Let $f := [0; a_1, a_2, \ldots ]$ be in $\wh K \setminus (K\cup \Theta_\alpha)$. 
Let $h$ be the largest non-negative integer for which there are 
arbitrarily large integers $n$ such that the $h$ polynomials 
$a_n, a_{n+1}, \ldots , a_{n+h-1}$ are of degree at most $k$. 
If $h \ge k$, then $c_{\alpha} (f) \le q^{-2 k - 2}$. 
If $h = k-1$, then there exists a polynomial $b$ of degree $k$ and infinitely many 
integers $n$ such that 
$a_n a_{n+1} \ldots  a_{n+h-1} b$ is a factor of $W$ and $\deg a_{n+h} > k$. This implies that 
$c_{\alpha} (f) \le q^{-2 (k-1) - k - 1} = q^{- 3 k + 1} \le q^{-2k - 1}$. 
If $0 < h \le k-2$, then there exists a polynomial $b$ of degree $k$ and infinitely many 
integers $n$ such that 
$b a_n a_{n+1} \ldots  a_{n+h-1} b$ is a factor of $W$, 
$\deg a_{n-1} > k$, and $\deg a_{n+h} > k$. This implies that 
$c_{\alpha} (f) \le q^{-2 h - 2k} \le q^{- 2 k -2}$. 

So, we are left with the case where all but finitely many $a_n$'s are polynomials of 
degree at least $k+1$. Since there are infinitely many pairs $Y^k, Y^{k+1}$ 
in the sequence of partial quotients of $\beta$, we get that $c_{\beta} (f) \le q^{- 1 - 2k}$
with equality, for instance, for $f = [0 ; \overline{Y^{k+2}}]$. 
All this implies that $\Hw(\beta) = q^{-2k - 1}$. 

Consequently, and taken also Proposition 4.2 (3) into account, we have shown that the function 
$\Hw$ takes any value $q^{-m}$, with $m=2$ or $m \ge 4$. To conclude, let $W''$ be a finite word 
over ${\cal A}_{\le 2}$ of even length, such that every polynomial of degree at most $2$ occurs 
in $W''$ and any two consecutive polynomials are of different degree. 
Let $\gamma$ be in $\QI$ whose period is given by $W''$. Then, it is easy to check that 
$\Hw (\gamma) = q^{-3} = c_\gamma ([0 ; \overline{Y^3}])$. 
This completes the proof of the theorem. 
\cqfd

\bigskip

\noindent {\it Proof of Theorem 1.5.} 
Let $k \ge 2$ and $\ell \ge 1$ be integers. 
Consider a cyclic de Bruijn word \cite{Bru46,Good46} of order $\ell$ over ${\cal A}_{=k}$, 
that is, a word $W$ of length $(\Card {\cal A}_{=k})^\ell$ such  that 
every word of length $\ell$ occurs exactly 
once in the prefix of length $(\Card {\cal A}_{=k})^\ell + \ell - 1$ of $WW$. 
Let $\alpha$ be in $\QI$ whose period is given by $W$. 

Let $j$ be a non-negative integer. Let $d, d'$ be positive integers different from $k$. 
Consider a factor $W_j = w_1 \ldots w_j$ of $W W$ of length $j$. 
It follows from Corollary 2.4 and Lemma 3.1 that 
$$
c_{\alpha} ( [0; \overline{Y^d, w_1, \ldots , w_j, Y^{d'}}] ) = q^{- 2 k j - \min\{d, k\} - \min\{d', k\}}. 
$$
Suitable choices of $d$ and $d'$ show that the function $c_\alpha$ can take 
every value between $q^{- 2 k j - 2}$ and $q^{- 2 k j - 2 k}$. 

It then remains for us to see which values of the form $q^{- 2 k j - 1}$ can be taken 
by the function $c_\alpha$.

Let $f := [0; a_1, a_2, \ldots ]$ be in $\wh K \setminus K$. 
Let $h$ be the largest non-negative integer for which there are 
arbitrarily large integers $n$ such that the polynomials 
$a_n, a_{n+1}, \ldots , a_{n+h-1}$ are of degree $k$. If $h \le \ell$, then 
one gets
$$
q^{- 2 k h - 2 k}  \le c_{\alpha} (f) \le q^{-2 k h - 2}. 
$$
Otherwise, we have
$$
c_{\alpha} (f) \le q^{-2 k h - k - 1}.
$$
This shows that the points $q^{- 2 k  - 1}, q^{- 4 k - 1}, \ldots , q^{-2 \ell k - 1}$ 
are not in the spectrum of $\alpha$. It remains for us to establish that, if $k$ is sufficiently large, then 
$q^{-2 h k - 1}$ is in $\Sp(\alpha)$ for every $h \ge \ell + 1$.

Let $Z$ be a word of length $\ell + 1$ over ${\cal A}_{=k}$ which is a factor of $WW$.
Let $P_1$ denote its last letter and write $Z = Z' P_1$. Let $P_2$ be a polynomial of degree $k$ such that 
$\deg (P_1 - P_2) \le k-1$ and 
with the property that neither 
$Z' P_2$, nor its mirror image, nor any of their twists by an element of $\kk^\times$ as described 
in Section 3, is a factor of $WW$. The existence of $Z, P_1, P_2$ is guaranted if $k$ is 
sufficiently large in terms of $\ell$ and the cardinality $q$ of $\kk$. 

If $f$ is a quadratic power series whose period is composed of $Y^d$
followed by the letters of $Z$, then one gets 
$$
c_{\alpha} (f) = q^{- 2 k \ell - \min\{d, k\} - 2 k + \deg (P_1 - P_2) }. 
$$
Choosing $d$ such that $ \min\{d, k\} = \deg (P_1 - P_2) +1$, 
this shows that there exists a power series $g$
such that $c_\alpha (g)$ is equal to $q^{-2 (\ell + 1) k - 1}$. 

A similar argument shows that $c_\alpha$ takes every value of the form 
$q^{-2 h k - 1}$, with $h \ge \ell + 2$. We omit the details. 
\cqfd

\medskip

The following example shows that the conjecture formulated by Parkkonen and Paulin, 
which is reproduced at the end of Section 1, does not hold. 
Set $N = q^2 - q$ and let $P_1, \ldots , P_N$ be all the (distinct) polynomials of degree $1$. 
Let $Q_1, \ldots , Q_N$ be distinct polynomials of degree $4$ and set 
$b_{2j - 1} = P_j, b_{2j} = Q_j$ for $j = 1, \ldots , N$. Let $\alpha$ be a quadratic power series 
whose period is given by $b_1, \ldots , b_{2N}$. 
Let $P$ be a polynomial and set $f_P := [0 ; \overline{P}]$. An easy computation shows that 
$$
\eqalign{
c_{\alpha} (f_P) & = q^{-5}, \quad \hbox{if $\deg P \ge 5$},  \cr
c_{\alpha} (f_P) & \le q^{-5}, \quad \hbox{if $\deg P = 4$},  \cr
c_{\alpha} (f_P) & = q^{-4}, \quad \hbox{if $\deg P = 3 $},  \cr
c_{\alpha} (f_P) & = q^{-3}, \quad \hbox{if $\deg P = 2$},  \cr
c_{\alpha} (f_P) & = q^{-4}, \quad \hbox{if $\deg P = 1$}.  \cr }
$$
Furthermore, we have 
$$
\Hw(\alpha) = q^{-3}.
$$
To see this, it is sufficient to observe that $c_{\alpha} (f) \le q^{-4}$ if $\alpha$ has infinitely many 
partial quotients of degree $1$ or $4$ and that $c_{\alpha} (f) \le q^{- 3}$ otherwise. 
 
It is tempting to conjecture that $\Hw(\alpha)$ is always attained at some 
power series $f_P := [0 ; \overline{P}]$, where $P$ is of minimal degree such that 
$P \not= b_i$ for every $b_i$ being in the periodic part of $\alpha$. 
However, such a conjecture also does not hold. 

\vskip 7mm

\noi {\bf Acknowledgements:}
I am very pleased to thank 
Fr\'ed\'eric Paulin 
for fruitful correspondence.

\vskip 23mm

\bigskip

\goodbreak

\centerline{\bf References}

\vskip 5mm

\beginthebibliography{999}

\bibitem{BPP}
A. Broise-Alamichel, J. Parkkonen, F. Paulin, 
Equidistribution and counting under equilibrium states in negatively curved spaces 
and graphs of groups. Applications to non-Archimedean Diophantine approximation. 
Book preprint (318 pages).


\bibitem{Bru46}
N. G. de Bruijn,  
{\it A combinatorial problem}, 
Indagationes Math. 8 (1946), 461--467.

\bibitem{Bu12}
Y. Bugeaud,
{\it Continued fractions with low complexity: 
Transcendence measures and quadratic approximation},
Compos. Math. 148 (2012), 718--750.

\bibitem{Bu14}
Y. Bugeaud,
{\it On the quadratic Lagrange spectrum},
Math. Z. 276 (2014), 985--999. 

\bibitem{CuFl}
T. W. Cusick and M. E. Flahive, 
The Markoff and Lagrange Spectra. 
Mathematical Surveys and Monographs, 
vol. 30, American Mathematical Society, Providence, RI, 1989.

\bibitem{CuMoPo96}
T. W. Cusick, W. Moran, and A. D. Pollington, 
{\it Hall's ray in inhomogeneous Diophantine approximation}, 
J. Austral. Math. Soc. Ser. A 60 (1996), 42--50. 

\bibitem{Fr75}
G. A. Freiman, 
Diophantine approximation and geometry of numbers (the Markov spectrum), 
Kalininski{\u\i }  Gos. Univ., Moscow, 1975. 

\bibitem{Gay17} 
D. Gayfulin, 
{\it Attainable numbers and the Lagrange spectrum}, 
Acta Arith. 179 (2017), 185--199. 

\bibitem{Good46}
I. J. Good, 
{\it Normal recurring decimals}, 
J. London Math. Soc. 21 (1946), 167--169. 

\bibitem{Hall47}
M. Hall,
{\it On the sum and product of continued fractions},
Ann. of Math. 48 (1947), 966--993.

\bibitem{HW} 
G. H. Hardy and E. M. Wright.
An introduction to the theory of numbers, 5th. edition, Clarendon Press, 1979.

\bibitem{HePa10}
S. Hersonsky and F. Paulin, 
{\it On the almost sure spiraling of geodesics in negatively curved manifolds}, 
J. Differential Geom. 85 (2010), 271--314.

\bibitem{HuMaUl15}
P. Hubert, L. Marchese, and C. Ulcigrai. 
{\it Lagrange Spectra in Teichmüller Dynamics via renormalization}, 
Geom. Funct. Anal. 25 (2015), 180--225. 

\bibitem{Las00}
A. Lasjaunias,
{\it A survey of Diophantine approximation in fields of power series},  
Monatsh. Math. 130 (2000), 211--229.

\bibitem{Lin18}
X. Lin,
{\it Quadratic Lagrange spectrum: I}, 
Math. Z. 
To appear. 

\bibitem{Math70}
B. de Mathan,
Approximations diophantiennes dans un corps local. 
Bull. Soc. Math. France Suppl. M\'em. 21 (1970), 93 pp. 

\bibitem{Mau03} 
F. Maucourant,
{\it Sur les spectres de Lagrange et de Markoff des corps imaginaires quadratiques}, 
Ergodic Theory Dynam. Systems 23 (2003), 193--205.

\bibitem{Mo18}
C. G. Moreira, 
{\it Geometric properties of the Markov and Lagrange spectra}. 
Preprint. {\tt https://arxiv.org/abs/1612.05782 }

\bibitem{PaPa11}
J. Parkkonen and F. Paulin,
{\it Spiraling spectra of geodesic lines in negatively curved manifolds},
Math. Z. 268 (2011), 101--142. Erratum: Math. Z. 276 (2014), 1215--1216. 


\bibitem{PaPa18}
J. Parkkonen and F. Paulin,
{\it On the nonarchimedean quadratic Lagrange spectra}.
Preprint. {\tt https://arxiv.org/abs/1801.08046}

\bibitem{Pau02}
F. Paulin,
{\it Groupe modulaire, fractions continues et approximation diophantienne en caract\'eristique $p$}, 
Geom. Dedicata 95 (2002) 65-85.

\bibitem{Pej16}
T. Pejkovi\'c, 
{\it Quadratic Lagrange spectrum}, 
Math. Z. 283 (2016), 861--869.

\bibitem{Per}
O. Perron,
Die Lehre von den Ketterbr\"uchen.
Teubner, Leipzig, 1929.

\bibitem{PiWo01}
C. G. Pinner and D. Wolczuk, 
{\it On the inhomogeneous Hall's ray of period-one quadratics}, 
Experiment. Math. 10 (2001), 487--495.

\bibitem{Schm00} 
W. M. Schmidt,
{\it On continued fractions and Diophantine approximation in
power series fields},
Acta Arith. 95 (2000), 139--166.

\endthebibliography

\vskip1cm

\noindent Yann Bugeaud  

\noindent  Institut de Recherche Math\'ematique Avanc\'ee, U.M.R. 7501 

\noindent  Universit\'e de Strasbourg et C.N.R.S.

\noindent 7, rue Ren\'e Descartes      

\noindent 67084 STRASBOURG  (FRANCE)

\vskip2mm

\noindent {\tt bugeaud@math.unistra.fr}

\bye